\documentclass[12pt]{amsart}

\usepackage{verbatim}
\usepackage{amsmath}
\usepackage{enumerate}
\usepackage{amssymb}
\usepackage{color}
\numberwithin{equation}{section}

\newcommand{\s}{\sqrt}

\setbox0=\hbox{$+$}
\newdimen\plusheight
\plusheight=\ht 0
\newcommand\+{\;\lower\plusheight\hbox{$+$}\;}

\setbox0=\hbox{$-$}
\newdimen\minusheight
\minusheight=\ht0
\renewcommand\-{\;\lower\minusheight\hbox{$-$}\;}

\setbox0=\hbox{$\cdots$}
\newdimen\cdotsheight
\cdotsheight=\plusheight
\newcommand\cds{\lower\cdotsheight\hbox{$\cdots$}}

\renewcommand{\s}{\text{sn}}
\renewcommand{\c}{\text{cn}}
\renewcommand{\d}{\text{dn}}

\begin{document}

\title[Representations of integers as sums of squares]%
{Recent progress in the study of representations of
integers as sums of squares }

\author{Heng Huat Chan and Christian Krattenthaler}

\address{
Department of Mathematics,
National University of Singapore,
2~Science Drive~2,
Singapore 117543}
\email{chanhh@math.nus.edu.sg}

\address{
Institut Girard Desargues,
Universit{\'e} Claude Bernard Lyon-I,
21, avenue Claude Bernard,
F-69622 Villeurbanne Cedex, France}
\email{kratt@euler.univ-lyon1.fr}

\begin{abstract}{In this article, we collect the recent results concerning the
representations of integers as sums of an even number of squares
that are inspired by conjectures of Kac and Wakimoto.
We start with a sketch of Milne's proof of two of these
conjectures. We also show an alternative route
to deduce these two conjectures
from Milne's determinant formulas for sums of $4s^2$,
respectively $4s(s+1)$, triangular
numbers. This approach is inspired by Zagier's proof
of the Kac--Wakimoto formulas via modular forms.
We end the survey with recent conjectures of the first author
and Chua.}
\end{abstract}

 \maketitle

 \section{Introduction}

The problem of finding explicit formulas for the number of
representations of an integer $n$ as a sum of $s$ squares is an
old one. The first formula of this kind is due to Legendre and
Gau\ss. If $r_s(n)$ denotes the number of representations of $n$ as
a sum of $s$ squares, Legendre and Gau{\ss} proved that
\begin{equation} \label{2squares} r_2(n)=4(d_1(n)-d_3(n)),\end{equation}
where $d_j(n)$ denotes the number of divisors of $n$ of the form
$4k+j$. For example, if $n$ is a prime $p$ of the form $4k+1$,
then $r_2(p) = 8$ since $d_1(p)=2$ and $d_3(p)=0$. On the other
hand, if $n$ is a prime $p$ of the form $4k+3$, then $r_2(p)=0$
since $d_1(p)=d_3(p)=1$. This, of course, leads to the well known
result of Fermat, which states that a prime $p$ is of the form
$x^2+y^2$ if and only if $p$ is of the form $4k+1$. Fermat's
result led mathematicians to explore and characterize primes of
the form $x^2+ny^2$, $n\geq 1$. For more information on such
characterizations, the reader is encouraged to consult the
excellent book by D.~A.~Cox \cite{Cox}.

Let $$\varphi(q)=\sum_{k=-\infty}^\infty q^{k^2}.$$
It is clear that
$$\varphi^s(q)=\sum_{k\geq 0} r_s(k)q^k.$$
As a result, to obtain expressions for $r_s(n)$, it suffices to
obtain expressions for $\varphi^s(q)$.
The first identity of this kind is due to Jacobi, namely,
\begin{equation}\label{2squares2}
\varphi^2(q) = 1+4\sum_{k=1}^\infty
(-1)^k\frac{q^{2k-1}}{1-q^{2k-1}}.\end{equation}
Note that \eqref{2squares} is a direct consequence of
 \eqref{2squares2}.
Using the theory of elliptic functions, Jacobi also found formulas for
$r_4(n)$, $r_6(n)$ and $r_8(n)$, namely,
\begin{align}
\varphi^4(q) &= 1+8\sum_{k=1}^\infty \frac{kq^k}{1+(-q)^k},
\label{4squares}\\
\varphi^6(q) &=1+16\sum_{k=1}^\infty
\frac{k^2q^k}{1+q^{2k}}-4\sum_{k=1}^\infty
(-1)^k\frac{(2k-1)^2q^{2k-1}}{1-q^{2k-1}},\\
\intertext{and}
\varphi^8(q) &= 1+16\sum_{k=1}^\infty
\frac{k^3q^k}{1-(-q)^k}.\label{8squares}\end{align}
From \eqref{4squares}, we find that
$$r_4(n) = 8\sum_{\substack{d|n \\
d\not\equiv 0 \ (\text{mod } 4)}} d,$$
and this immediately implies that every positive integer is
a sum of four squares, a famous result of Lagrange.

We call series of the type $$A+B\sum_{k\geq
1}a_k\frac{q^k}{1-q^k}$$ {\it generalized Lambert series}. Note
that for even $s\leq 8$, we are able to express $\varphi^s(q)$ in
terms of generalized Lambert series. This does not seem possible
when $s=10$. In fact, Liouville showed that \begin{align*}
\varphi^{10}(q) &= 1+\frac{4}{5}\sum_{k=1}^\infty
(-1)^{k-1}\frac{(2k-1)^4q^{2k-1}}{1-q^{2k-1}}+\frac{64}{5}\sum_{k=1}^\infty
\frac{k^4q^k}{1+q^{2k}} \\ &\qquad +
\frac{32}{5}q\varphi^2(q)\varphi^4(-q)\psi^4(q^2),\end{align*}
where $$\psi(q) = \sum_{k=0}^\infty q^{k(k+1)/2}.$$ Indeed for any
even $s\geq 10$, $\varphi^s(q)$ is a sum of generalized Lambert
series and a ``cusp form''.

Recently, new formulas for $r_s(n)$ were discovered. One common
feature of these formulas is the absence of ``cusp
forms''.\footnote{Note that the coefficients of $q^n$ in Lambert
series can be calculated once we know the factorization of $n$. In
general, this is impossible for cusp forms. Hence, these new
formulas are more ``effective'' if one wants to determine
$r_s(n)$.}  The new formulas involve only generalized Lambert
series. The purpose of this article is to describe these recent
discoveries.

Before we proceed with our discussion, we make the following
observation: It is known that \cite[p. 43, Entry 27 (ii)]{PartIII}
\begin{equation}\label{ts} 4e^{\pi i/(2\tau)}
\psi^2(e^{-2\pi i/\tau}) =
\frac{\tau}{i}\varphi^2(-e^{\pi i\tau}).
\end{equation}
Suppose we
have a relation
$$4^sq^{s/2}\psi^{2s}(q^2) = F(L_1(q^2),L_2(q^2),\dots, L_m(q^2)),\quad \text{with }
q=e^{\pi i \tau},$$ where each $L_j$ is a generalized Lambert
series or a product of generalized Lambert series satisfying
$$L_j(e^{-2\pi i/\tau}) = \left(\frac{\tau}{i}\right)^s L_j^*(-e^{\pi i\tau})$$
 for some $L_j^*(-q)$ (which is also a generalized Lambert series or a
 product of generalized Lambert series), then
we would have
 $$\varphi^{2s}(-q) = F(L_1^*(-q),\dots, L_m^*(-q)).$$
 Conversely if we have a formula for sums of squares,
 we will have a formula for sums of triangular numbers.
We illustrate the above observation by the following identities:

Suppose for $q=e^{\pi i \tau}$, we have
\begin{align} \label{8t} 4^4e^{2\pi i\tau}\psi^8(e^{2\pi i \tau}) &=16\sum_{k\geq
0}\frac{k^{3}q^{2k}}{1-q^{4k}}\\
&=\frac{16}{15}\left(E_4(\tau)-E_4(2\tau)\right),\notag\end{align}
where \begin{equation}\label{tE4} E_4(\tau) = 1+240\sum_{k\geq
1}\frac{k^3e^{2\pi i k\tau}}{1-e^{2\pi i k\tau}}.\end{equation}
The Eisenstein series $E_4(\tau)$ satisfies the transformation
formula \cite[p. 24, Ex. 12]{Apostol}
$$E_4\left(-\frac{1}{\tau}\right) =\tau^4 E_4(\tau).$$
Replacing $\tau$ by $-1/\tau$ in \eqref{8t}, the left hand side of
\eqref{8t} is $\tau^4\varphi^8(-q)$ by \eqref{ts}, with $q=e^{\pi
i \tau}$. Now, by \eqref{tE4}, we have
\begin{align*}
\frac{16}{15}&\left(E_4(-1/\tau)-E_4(-2/\tau)\right)=
\tau^4\frac{16}{15}\left(E_4(\tau)-\frac{1}{2^4}E_4(\tau/2)\right)\\
&=\tau^4\left(1+128\sum_{k\geq
1}\frac{k^3q^k}{1-q^k}-16\sum_{k\geq
1}\frac{(2k+1)^3q^{2k+1}}{1-q^{2k+1}}\right).\end{align*} Hence,
we conclude that
$$\varphi^8(-q)=1+16\sum_{k=1}^\infty \frac{k^3(-q)^k}{1-q^k}.$$
Replacing $q$ by $-q$, we obtain the formula for sums of 8
squares. For more details of the relations between formulas
associated with squares and triangular numbers see \cite{Liu} and
\cite{OnoJNT}.

Note that the identity for $\psi^{8}(q)$ is much simpler than
that for $\varphi^{8}(q)$. This is in fact a general
phenomenon: the identity for $\psi^{2s}(q)$ will be much simpler than
that for $\varphi^{2s}(q)$ for
 any $s\in\mathbb N$. For the rest of this article,
we will therefore only present identities associated with $\psi(q)$.

\section{The formulas of Kac and Wakimoto}

In 1994, V.~G.~Kac and M.~Wakimoto \cite{KacWakimoto} conjectured
that
\begin{equation} \label{KW} t_{4s^2}(n) =\frac{1}{s!}\frac{4^{-s(s-1)}}{\prod_{j=1}^{2s-1}
j!}\sum_{\substack{a_1,\cdots, a_s\in \mathbb N,\ \text{$a_i$ odd} \\
r_1,\cdots , r_s\in \mathbb N,\ \text{$r_i$ odd} \\
a_1r_1+\cdots a_sr_s=2n+s^2}} a_1\cdots
a_s\prod_{i<j}(a_i^2-a_j^2)^2\end{equation}
and
\begin{equation} \label{KW2} t_{4s(s+1)}(n) =\frac{1}{s!}
\frac{2^s}{\prod_{j=1}^{2s}
j!}\sum_{\substack{a_1,\cdots, a_s\in \mathbb N \\
r_1,\cdots , r_s\in \mathbb N,\ \text{$r_i$ odd} \\
a_1r_1+\cdots a_sr_s=n+\frac{1}{2}s(s+1)}} (a_1\cdots
a_s)^3\prod_{i<j}(a_i^2-a_j^2)^2.\end{equation}
These formulas follow from a conjectural affine denominator formula
for simple Lie superalgebras 
 of type $Q(m)$.
(For the definition of $Q(m)$, see \cite{KacSuper}.)

Identities \eqref{KW} and \eqref{KW2} were first proved by S.~C.~Milne
\cite{Milne}, using
results on continued fractions and elliptic functions. For example,
Milne showed using Schur functions, that \eqref{KW} is a consequence of his
determinant formula \cite[(5.107)]{Milne}
\begin{equation}\label{Milne}
\left(q\psi^4(q^2)\right)^{s^2}=\frac{4^{-s(s-1)}}{\prod_{j=1}^{2s-1} j!}
\text{det}(C_{2(u+v-1)-1})_{1\leq u,v
\leq s}\ ,\end{equation}
where
$$C_{2j-1} = \sum_{r=1}^\infty
\frac{(2r-1)^{2j-1}q^{2r-1}}{1-q^{2(2r-1)}}, \quad j\geq 1.$$
We now briefly describe Milne's proof of \eqref{Milne}.

Milne first showed that if $\s(u):=\s(u,\mathbf{k})$,
$\d(u):=\d(u,\mathbf{k})$, and
$\c(u):=\c(u,\mathbf{k})$ are the classical
Jacobi elliptic functions, 
then \cite[(2.44), (2.68)]{Milne}
\begin{align}\label{fourier}
\frac{\s(u)\c(u)}{\d(u)} &=\frac{1}{\mathbf{k}^2}\sum_{m\geq 1}
\frac{2^{2m+2}(-1)^{m-1}}{z^{2m}} C_{2m-1}
\frac{u^{2m-1}}{(2m-1)!}\\
&=:\sum_{m\geq 1}c_m \frac{u^{2m-1}}{(2m-1)!},
\notag
\end{align}
where \begin{equation}\label{zx} z=\varphi^2(q)
\quad\text{and}\quad
\mathbf{k}^2=16q\frac{\psi^4(q^2)}{\varphi^4(q)}.\end{equation}
Milne then showed that
\begin{multline}\label{sc/d} \int_0^\infty
\frac{\s(u)\c(u)}{\d(u)}e^{-u/t} \,du \\=
\frac{t^2}{1+(4-2\mathbf{k}^2)t^2+
{\raise-5pt\hbox{\Huge $\mathbf{K}$}_{n=2}^\infty}
\dfrac{-(2n-1)(2n-2)^2(2n-3){\mathbf{k}}^4t^4}{1+(2n-1)^2
(4-2\mathbf{k}^2)t^2}}.
\end{multline}
Here, $\mathbf K_{n=2}^{\infty}$ is the notation for continued fractions,
$${\raise-5pt\hbox{\Huge $\mathbf{K}$}_{n=2}^\infty}
\frac {a_n} {b_n}:=
\cfrac {a_2}{b_2+
\cfrac {a_3}{b_3+
\cfrac {a_4}{b_4+\raise-5pt\hbox{$\ddots$}
}}}\ .
$$
Using \cite[Theorem 3.4]{Milne} and \eqref{sc/d}, Milne deduced
the Hankel determinant evaluation \cite[(4.9)]{Milne}
\begin{equation}\label{hankel}
H_n^{(1)}(\{c_m\}):=\text{det}\left(\begin{matrix} c_1 & c_2 &\cdots &c_n \\
c_2 &c_3 &\cdots &c_{n+1} \\
\vdots &\vdots &\ddots &\vdots \\
c_n &c_{n+1}&\cdots &c_{2n-1}\end{matrix}\right)
=({\mathbf k}^2)^{n(n-1)}
\prod_{r=1}^{2n-1}r!.\end{equation}
Simplifying the left hand side of \eqref{hankel} using
the definition \eqref{fourier} of the $c_i$'s and making use of
the relations \cite[(3.66), (5.11)]{Milne}
$$H_n^{(1)}(\{t^ma_m\})=t^{n^2}H_n^{(1)}(\{a_m\}),$$
and \eqref{zx},
we deduce \eqref{Milne}.


%
We now describe a simplification of Milne's Schur function argument
that allowed him to deduce \eqref{KW} from \eqref{Milne}.
As a side result, we also obtain a new expression for $t_{4s^2}(n)$
(see \eqref{CC} below).

In a recent paper \cite{Zagier}, D.~Zagier gave a direct proof of
the above formulas of Kac and Wakimoto using the theory of modular
forms. In that paper, he constructed a certain map sending the
monomials $X_1^{k_1-1}\cdots X_s^{k_s-1}$ (here, the $X_i$'s are
indeterminates) to the product of Eisenstein series
$g_{k_1}^+\cdots g_{k_s}^+$ (the quantities $g_{2j}^+$ being, up
to scaling, the quantities $C_{2j-1}$ in Milne's formula). It
turns out that if we apply a variant of that map, $\Phi_s$ say,
defined by sending the product $X_1^{2k_1-1}\cdots X_s^{2k_s-1}$
to the product $C_{2k_1-1}\cdots C_{2k_s-1}$, in Milne's formula,
then we get the new formula
\begin{multline} \label{CC}
t_{4s^2}(n) = \frac{(-1)^{s(s-1)/2}}{4^{s(s-1)}
\prod_{j=1}^{2s-1}j!}\\
\times
\sum_{\substack{\text{$a_i, r_i\in \mathbb N$ odd} \\
a_1r_1+\cdots +a_sr_s=2n+s^2}} a_1a_2^3\cdots
a_s^{2s-1}\prod_{1\leq i<j\leq s} (a_i^2-a_j^2).\end{multline}
This is seen as follows: By series expansion, we have
\begin{align*}
C_{2j-1}&= \sum_{r=1}^\infty
(2r-1)^{2j-1}q^{2r-1}
\sum _{k=0} ^{\infty}q^{2k(2r-1)}\\
&=
\sum _{r,k\ge1} ^{}(2r-1)^{2j-1}q^{(2k-1)(2r-1)}\\
&=
\underset {a\mid m}{\sum _{m\text{ odd}} ^{}}a^{2j-1}q^m.
\end{align*}
Thus, we obtain
\begin{align*}
\Phi_s(X_1^{2k_1-1}\cdots
X_s^{2k_s-1})&=C_{2k_1-1}\cdots C_{2k_s-1}\\
&=\underset {a_1\mid m_1,\
\dots,\ a_s\mid m_s}{
\sum _{m_1,\dots,m_s\text{ odd}}
^{}}q^{m_1+\dots+m_s}a_1^{2k_1-1}\cdots a_s^{2k_s-1}.
\end{align*}
Returning to Milne's formula \eqref{Milne}, this implies that
\begin{align*}
\det\big(&C_{2(u+v-1)-1}\big)_{1\leq u,v\leq s}=
\Phi_s\left(\det\left(X_u^{2(u+v-1)-1}\right)_{1\leq u,v\leq s}\right)\\
&=\Phi_s\left(
\prod _{i=1} ^{s}X_i^{2i-1}
\det\left(X_u^{2(v-1)}\right)_{1\leq u,v\leq s}\right)\\
&=\Phi_s\left((-1)^{s(s-1)/2}
\prod _{i=1} ^{s}X_i^{2i-1}
\prod_{1\leq i<j\leq s}(X_i^2-X_j^2)\right)\\
&=(-1)^{s(s-1)/2}\underset {a_1\mid m_1,\
\dots,\ a_s\mid m_s}{
\sum _{m_1,\dots,m_s\text{ odd}}
^{}}q^{m_1+\dots+m_s}\prod _{i=1} ^{s}a_i^{2i-1}
\prod_{1\leq i<j\leq s}(a_i^2-a_j^2),
\end{align*}
where we have used the Vandermonde determinant evaluation to evaluate
the determinant in going from the second to the third line. Now a
comparison of coefficients of $q^{2n+s^2}$ leads us to
\eqref{CC}.

We now show that \eqref{KW} follows from \eqref{CC} using an elementary
combinatorial argument.

For each positive integer $s$, let  $$P_s(X_1,\dots, X_s) =
\prod_{i=1}^s X_i\prod_{i<j} (X_i^2-X_j^2)^2$$ and
$$P_s'(X_1,\dots, X_s) = \prod_{i=1}^s
X_i^{2i-1}\prod_{i<j}(X_i^2-X_j^2).$$
For positive integers $s$ and $m$, let
$$R_s(m,P(X_1,\dots, X_s)) =\sum_{\substack{\text{$a_i, r_i\in \mathbb N$ odd} \\
a_1r_1+\cdots +a_sr_s=m}}P(a_1,\dots, a_s).$$ We want to show
that
$$R_s(m,P_s)=(-1)^{s(s-1)/2}s!\,R_s(m,P'_s).$$
Now, if $S_s$ denotes the symmetric group of $s$ elements and
$\sigma\in S_s$, then
\begin{equation*}
R_s(m,P'_s(X_{\sigma(1)},\dots,X_{\sigma(s)}))
=R_s(m,P'_s(X_1,\dots , X_s)),\end{equation*}
since, whenever $(a_1,\dots,a_s)$ is an $s$-tuple of odd nonnegative integers
for which there are
$r_1,\dots r_s$ such that $a_1r_1+\dots+a_sr_s=m$,
then $(a_{\sigma(1)},\dots,
a_{\sigma(s)})$ is an $s$-tuple with the same property.
Hence,
$$R_s\bigg(m,\sum_{\sigma\in S_s}
P'_s(X_{\sigma(1)},\dots,X_{\sigma(s)})\bigg) =
s!\,R_s(m,P'_s(X_1,\dots , X_s)).$$ On the other hand,
\begin{align*}
R_s\bigg(m,\sum_{\sigma\in S_s}&
P'_s(X_{\sigma(1)},\dots,X_{\sigma(s)})\bigg)\\
&= R_s\bigg(m,
\text{det}(X_j^{2i-2})_{1\le i,j\le s}
\prod _{i=1} ^{s}X_i\prod_{1\le i<j\le s}(X_i^2-X_j^2)\bigg) \\
&= (-1)^{s(s-1)/2}R_s(m,P_s(X_1,\dots, X_s)),\end{align*}
by the Vandermonde determinant evaluation.
Thus, we have shown how \eqref{KW} follows from
Milne's determinant formula \eqref{Milne} by passing via
\eqref{CC}, thereby providing an alternative to Milne's
(somewhat more involved) Schur function argument.

\section{A conjecture for the sum of $8s$ triangular numbers}

We now state Milne's determinant formula for $4s(s+1)$ triangles:
\begin{equation}\label{Milne2} \left(16q\psi^4(q^2)\right)^{s(s+1)}
=\left(2^{s(4s+5)}\right)\prod_{j=1}^{2s} (j!)^{-1}
\text{det}\left(D_{2(u+v-1)+1}\right)_{1\leq u,v\leq s},
\end{equation}
 where
$$D_{2j+1} = \sum_{r=1}^\infty \frac{r^{2j+1}q^{2r}}{1-q^{4r}},
j\geq 1.$$
This formula led to the first proof of \eqref{KW2}. Using
arguments analogous to the ones
given in the last section, one can deduce \eqref{KW2} from
\eqref{Milne2}.

When $s=2$, this leads to the following beautiful formula:
$$q^6\psi^{24}(q^2) = \frac{1}{72}\left(T_8T_4-T_6^2\right),$$
where
$$T_{2k}(q) :=\sum_{n=1}^\infty \frac{n^{2k-1}q^{2n}}{1-q^{4n}},\ k>1.$$
 Note the resemblance of this formula with the  well-known formula
$$q\prod_{n=1}^\infty (1-q^n)^{24} =
\frac{1}{1728}\left(E_4E_8-E_6^2\right),$$ where the $E_i$'s are the
classical Eisenstein series.\footnote{This was probably first observed
by F.~G.~Garvan.}
Note that, as indicated in the introduction, one obtains Milne's new
formula for 24
squares, namely, if
\begin{align*} S_4(q) &=
1+16\sum_{k=1}^\infty\frac{k^3q^k}{1-(-q)^k},\\
S_6(q) &= 1-8\sum_{k=1}^\infty \frac{k^5q^k}{1-(-q)^k}, \\
\intertext{and}
 S_8(q) &= 17+32\sum_{k=1}^\infty
\frac{k^7q^k}{1-(-q)^k},
\end{align*}
then \cite[Theorem 1.6, (1.25)]{Milne}
\begin{equation*}
\varphi^{24}(q) =
\frac{1}{9}\left\{S_4(q)S_8(q)-8S_6^2(q)\right\}.
\end{equation*}
Comparing this with the ``old'' formula
\begin{align*}\varphi^{24}(q) &= 1+\frac{16}{691}E_{11}(q)
+\frac{33152}{691}qf^{24}(q)-\frac{65536}{691}q^2f^{24}(-q^2),
\end{align*}
where $$f(-q) = \prod_{k=1}^\infty (1-q^k),$$ we find that Milne's
formula requires less terms. Moreover, if we know the
factorization of $n$, then we can calculate $r_{24}(n)$ explicitly
from the new formula, since the terms are all Eisenstein series.
The new formula for $24$ squares and a recent paper of Z.-G. Liu
\cite{Liu} led the first author and Chua \cite{32squares} to
formulate the following conjecture for $8s$ triangular numbers:
\medskip

\noindent{\bf Conjecture.}
For any positive integer $s>1$, we have
$$q^{2s}\psi^{8s}(q^2) = \sum_{\substack{m+n=2s \\ m\ge n\geq 2}}
a_{m,n}T_{2m}T_{2n},$$ for some rational numbers $a_{m,n}.$
\medskip

For a fixed $s$ one can verify the corresponding
identity. For example, when $s=4$ we are led to the following new
identity:
\begin{equation}\label{32} q^8\psi^{32}(q^2) = \frac{1}{75600}\left(\frac{25}{4}T_{10}(q)T_{6}(q)-
\frac{21}{4}T_8^2(q)-T_4(q)T_{12}(q)\right).\end{equation}
Note that the above identity does not follow from any of the formulas
of Kac and Wakimoto or Milne, since 32 is not of the form $4s^2$ or $4s(s+1)$.

The first proof of this result proceeds by expressing the
$T_{2m}$'s in terms of $\mathbf k^2$ and $z$
(see \eqref{zx} for their definitions).
The corresponding expressions are found by using
the theory of modular forms, as well as the following recurrence
satisfied by $T_{2m}$'s:
\begin{equation*}
T_{2n+8}(q) = T_2(q)T_{2n+6}(q)
+12\sum_{j=0}^n \left(\begin{matrix} 2n+4 \\
2j+2\end{matrix}\right) T_{2j+4}(q) T_{2n-2j+4}(q),\end{equation*}
where
 $$T_2(q) = 1+24\sum_{j=0}^\infty \frac{jq^{2j}}{1+q^{2j}}.$$
The above recurrence follows from
the differential equation
satisfied by the Jacobi elliptic function $M:=M(u)=\text{sn}^2(u)$,
namely,
$$\left(\frac{dM}{du}\right)^2 = 4M(1-M)(1-{\mathbf k}^2M).$$

We end this article with a sketch of a new proof of  \eqref{32}.
It is known that the constant term of $M(u)$ is $T_{4}$ and that $$T_4=
q^2\psi^8(q^2).$$ One can verify that
\begin{align}\label{32t} 3840M^4 = &(7M^{(2)}+68T_4)(M^{(2)}-4T_4)\\
&\quad -9M^{(3)}M^{(1)}+M(2M^{(4)}+128T_6), \notag\end{align}
where $M^{(i)}$ is the $i$th derivative of $M$ with respect to $u$.

By comparing the constant term on both sides we are immediately led to
\eqref{32}. The above identity is motivated by recent work of
the first author and
Liu \cite{quintic}, where a proof of an identity  similar to
\eqref{32t} is illustrated.

\noindent{\it Acknowledgements.}  It is our pleasure to thank Professors Bruce C.
Berndt and S. Milne for pointing out several misprints in the preliminary
version of this article and for their fruitful suggestions.

%


\begin {thebibliography}{9}

\bibitem{Apostol} T.M. Apostol, {\it Modular functions and
Dirichlet series in Number Theory}, Springer-Verlag, 2nd ed., New
York, 1990.

\bibitem{PartIII} B.C. Berndt, {\it Ramanujan's Notebooks, Part
III}, Springer-Verlag, New York, 1991.

\bibitem{32squares} H. H. Chan and K.S. Chua, {\it Representations of integers
as sums of 32 squares}, The Ramanujan J. {\bf 7} (2003), 79--89.

\bibitem{quintic} H. H. Chan and Z.-G. Liu, {\it Elliptic functions to the
quintic base,} preprint.

\bibitem{Cox} D. A. Cox,
{\it Primes of the form $x^2+ny^2$}, John
Wiley \& Sons, New York, 1989.

\bibitem{KacSuper} V. G. Kac, Lie superalgebras, Adv.\ Math.\
{\bf 26} (1997), pp. 8--96.
\bibitem{KacWakimoto} V. G. Kac and M. Wakimoto, {\it Integrable highest
weight modules over affine superalgebras and number theory,} in: {\it Lie
Theory and Geometry,} in honor of Bertram Kostant
(J.~L.~Brylinski, R.~Brylinski, V.~Guillemin
and V.~Kac, eds.) {\bf 123}, Prog.\ in Math., Birkh\"auser Boston,
Inc., Boston, MA., 1994, pp. 415--456.

\bibitem{Liu} Z.-G. Liu, {\it On the representations of integers as sums of
squares,} in {\it $q$-series with applications to Combinatorics, Number Theory,
and Physics } (B.~C.~Berndt and K.~Ono, eds.), {\bf 291}, Contemp.\ Math.,
AMS, Providence, R.I., 2001, pp.~163--176.

\bibitem{Milne} S. Milne, {\it Infinite families of exact sums of squares
formulas, Jacobi elliptic functions, continued fractions, and Schur
functions}, The Ramanujan J. {\bf 6} (2002), 7--149.

\bibitem{OnoJNT} K. Ono, {\it Representations of integers as sums of
squares,} J. Number Theory {\bf 95} (2002), 253--258.


\bibitem{Zagier} D. Zagier, {\it A proof of the Kac-Wakimoto affine
denominator formula for the strange series,} Math.\ Res.\ Lett.\
{\bf 7} (2000), 597--604.

\end{thebibliography}
\end{document}